\documentclass{amsart}
\newtheorem*{thm}{\rm THEOREM}
 \newtheorem{lem}{\rm LEMMA}
\newtheorem{cond}{Condition}

 \theoremstyle{remark}
\newtheorem{rem}{REMARK}

 \numberwithin{equation}{section}

\begin{document}

\title{On the zeros of functions in the Selberg class}
\author{Anirban Mukhopadhyay, Kotyada Srinivas, Krishnan Rajkumar}
\address{Institute of Mathematical Sciences,
CIT Campus, Tharamani, Chennai 600 113, India}
\email[Anirban Mukhopadhyay]{anirban@imsc.res.in}
\email[Kotyada Srinivas]{srini@imsc.res.in}
\email[Krishnan Rajkumar]{rkrishnan@imsc.res.in}
\begin{abstract} It is proved that under some suitable conditions, the 
degree two functions in the Selberg class have infinitely many
zeros on the critical line.
 \end{abstract}

 \maketitle
\hfill \textit{To Professor Swadheen Pattanayak with deep regards}

\section{\emph{ Introduction}}
Let $\mathcal{S}$ denote the Selberg class of functions (see \cite{Per}).
Let $F \in \mathcal{S}$ be a function of {\it degree} $d > 0 $ (defined later). From the theory of integral functions of finite order,
 it is easy to show that $F$ has infinitely many zeros in the critical strip $0 \le \sigma \le 1$. In fact, a formula
 analogous to the Riemann - von Mangoldt formula for the Riemann zeta-function holds for $F$ (see \cite{Per}). It is therefore
natural to ask the following 

\medskip

\noindent \textbf{Question}: \textit{Does $F$ admit infinitely many zeros on the critical line $\sigma = \frac{1}{2}$}?

\medskip

The question has an affirmative answer for all degree $1$ functions and all known functions of degree $2$ in $\mathcal{S}$.

In \cite{K-P1},\cite{K-P2} Kaczorowski and Perelli show that the only functions of degree $0 < d < 5/3$  are the Riemann 
zeta-function $\zeta(s)$, the Dirichlet $L$-functions $L(s,\chi)$ associated with primitive characters $\chi\,$ and their
imaginary shifts $L(s+i\,\omega,\chi)$ with real $\omega$. (Recently
we learnt from Professor A. Ivi\'c in a personal communication that this result has been extended by the same authors upto $0 < d < 2$.)
In the case of $\zeta(s)$, Hardy \cite{Har} showed that it has an infinity of zeros on the critical line. A similar argument works for $L(s,\chi)$. In fact, for these functions it is known that a positive proportion of zeros lie on the critical line (see \cite{Sel}, \cite{lev}, \cite{con}, \cite{Zur}).

The known examples of degree $2$ functions are product of two degree 1 functions, for example $\zeta^2(s)$, $L(s,\chi_1)L(s,\chi_2)$; 
the $L$-function $L_f(s)$ associated with a suitable normalized modular form $f$ which is either holomorphic
or a Maass waveform ; the Dedekind zeta-function of quadratic
number field $K$, $\zeta_K(s)$; and the imaginary shifts of any of the above functions which are entire, for example $L_f(s)$, $L(s,\chi_1)L(s,\chi_2)$ when $\chi_1$ and $\chi_2$ are non-principal.
 In all these cases, it is known that they have infinitely many zeros on the critical line.
 In fact, it is known that a positive proportion of zeros lie on the critical line (see
 \cite{Haf1},\cite{Haf2}). However, since it is not yet known if this list is complete for
 degree 2, the question has not yet been anwered in general.

 %In this note we wish to show that under some conditions $F$ admits infinitely many 
%zeros on the {\it critical line.  } Such resuls are well-known for all the
%known degree two functions in the Selberg class. 
In this context, Gritsenko in \cite{Grit} starts with $F_1,F_2,\ldots F_N$ which are distinct primitive functions of degree $2$
in $\mathcal{S}$ which satisfy the following conditions:
\begin{cond}\label{cond1}
The following asymptotic formula holds as $ x \to \infty  $:
$$\sum_{p \le x}|a_j(p)|^2 \log p = x + O(x \log^{-1} x),$$
where $a_j(n)$ are the Dirichlet coefficients of $F_j(s), \quad j=1,2,\ldots,N.$
\end{cond}

\begin{cond}\label{cond2}
For any $j$, $1\le j \le N$, there exists a positive constant $c_{F_j}$ such that
$$|1+a_j(p)p^{-\frac{1}{2}-it}+a_j(p^2)p^{-1-2it}+\cdots|> c_{F_j},$$
where $t$ is an arbitrary number in the interval $ [T,2T+1], \ p$ is an arbitrary prime number and
$ T \geq T_0 > 0,  \ T_0 $ is a sufficiently large absolute constant.
\end{cond}

\begin{cond}\label{cond3}
The following asymptotic formula holds  as $ x \to \infty $ :
$$\sum_{n \le x}|a_j(n)|^2 n^{-1} = A_{F_j} \log x + B_{F_j} + O(\log^{-4} x),$$
where $A_{F_j} > 0$ and $B_{F_j}$ are constants depending  only on $F_j(s), \; j=1,2,\ldots,N.$
\end{cond}

He then assumes the truth of the following conjecture of Selberg:

\begin{cond}\label{cond4}
The following asymptotic formula holds as $ x \to \infty   $:
$$\sum_{p \le x}a(p)\overline{a'(p)}p^{-1} = O(1),$$
where $a(n)$ and $a'(n)$ are the Dirichlet coefficients of two \textbf{distinct} 
primitive functions $F$ and $F'$ in $\mathcal{S}$.
\end{cond}
He then shows that \textit{the function $\mathcal{F}$ given by
$$\mathcal{F}(t)=\sum_{j=1}^{N} b_j Z_{F_j}(t),$$
has at least $T \exp(\sqrt{\log \log \log T})$ zeros of odd order in $[T,2T]$, 
where $Z_{F_j}(t)$ is defined as in Section \ref{proof} and $b_j$ are real with $b_1\neq 0$.}

\medskip

An important result in this context is by Bombieri and Hejhal (cf. \cite{B-H}),
which deals with a class of functions $L(s)$ very close to the Selberg class. 
Roughly speaking, their result is as follows. \textit{If $L_1(s),...,L_N(s)$ satisfy the same
functional equation, are orthogonal in the sense of Selberg, satisfy the Riemann
Hypothesis and a certain weak form of the pair correlation conjecture, then almost
all zeros of the linear combination $$\sum_{j=1}^N b_j L_j(s),$$ where $b_j$ are real,
lie on the critical line $\sigma = 1/2$ and are simple.}
\medskip

Unconditionally, the question seems to be difficult for a general function in $\mathcal S$ of degree $d\ge 2$. An
attempt has been made in this note to answer the question for functions of degree $2$ under some mild conditions.
\medskip

To state the main theorem we need some preparation. We recall that every 
$ F \in \mathcal{S}$ satisfies a functional equation of the type
$$\Phi (s) = \omega \bar{\Phi} (1-s), 
$$
where
\begin{equation}\label{int0e}
\Phi(s) = Q^s \ \prod_{j=1}^{r} \Gamma(\lambda_js + \mu_j) F(s) 
=\gamma(s) F(s),
\end{equation}
say, with $ Q > 0,\lambda_j > 0, \Re \mu_j \geq 0 \ \  \rm{and} \ 
|\omega| = 1.$ Here
$ \bar{f}(s) = \overline{f(\bar{s})}.$ 

The functional equation can be rewriten as
$$F(s)=\Delta_F(s) \bar{F}(1-s),
$$
where
$$\Delta_F(s)=\omega \frac{\bar{\gamma}(1-s)}{\gamma(s)} \cdot
$$

 We also recall that $\log F(s)$ is a Dirichlet series with coefficients $b(n)$ 
satisfying $b(n) = 0$ unless $n = p^m,$ for some prime $p$ and integer 
$ m \geq 1$. This implies that $F$ has an Euler product expansion and hence the Dirichlet
 coefficients $a(n)$ of $F$ are multiplicative. Multiplicativity of the coefficients is
crucial in the proof of our theorem.

From now onwards the function $F\in \mathcal{S}$ will always
be denoted by $ F(s) = \sum_{n=1}^{\infty}\frac{a(n)}{n^s}$, $s=\sigma +it$ with 
the following {\it data}, which are invariants, associated with it (see \cite{Per}):
%$  \hbox{the degree} \ d :=  d_F=2 \sum_{j=1}^r \lambda_j , \  
% \hbox{the} \ \xi-invariant \xi:= \xi_F=2 \sum_{j=1}^r \left(\mu_j - 1/2 \right) $
%and $ \hbox{the conductor }  q := q_F= (2\pi)^{d} Q^2 \prod_{j=1}^r \lambda_j^{2 \lambda_j}.$
the degree $d = 2 \sum_{j=1}^r \lambda_j$, 
the $\xi$-invariant $\xi = 2 \sum_{j=1}^r \left(\mu_j - 1/2 \right)$
and the conductor $q = (2\pi)^{d} Q^2 \prod_{j=1}^r \lambda_j^{2 \lambda_j}$.
The implied constants in the $O$-terms appearing at various places may depend only on the conductor $q$.
\noindent
The aim of this note is to establish the following
\begin{thm}
Let $F \in S$ with $d=2$, $\xi$ real and $\sqrt q$ irrational 
and suppose the Dirichlet coefficients  $a(n)$ of $F$ satisfy
\begin{equation}\label{int1e}
\sum_{n \le x}|a(n)|^2 = O(x),
\end{equation}
then 
$F(s)$ has infinitely many zeros on the critical line $\Re s = \frac{1}{2}$.
\end{thm} 
\begin{rem}
The condition \eqref{int1e} can be derived from Condition \ref{cond3} by a straight forward application of partial summation.
\end{rem}
\begin{rem}
The presence of the factor $t^{ -i \Im \xi }$ in \eqref{asym1e}  makes it difficult to apply Lemma \ref{lem1} to evaluate
the integral in \eqref{pro4e}. Therefore, we imposed the condition that  $\xi$ is real.
Moreover, this condition is benign for functions $F$ which are regular at $s=1$, 
as a suitable imaginary shift of $F$ 
can be made to satisfy this and hence the conclusion of the theorem will hold for $F$ also.
\end{rem}
\begin{rem}
The conditions of the theorem are satisfied for the functions \\
\mbox{$L(s,\chi_1) L(s,\chi_2)$} where $\chi_1$ and $\chi_2$ are primitive Dirichlet characters to the modulus $q_1$ and $q_2$ respectively provided $q_1 q_2$ is not a perfect square. The conditions also hold for $L_f(s)$ for a
 cuspidal modular form $f$ of level $N$ if $N$ is not a perfect square.
\end{rem}

\begin{rem}
The conditions of the theorem are {\it not} satisfied for $\zeta^2(s)$ or $L^2(s,\chi)$. However, since our proof follows closely  Hardy - Littlewood proof ( as given in \cite{T}, pp 260--262 ) of Hardy's theorem for $\zeta(s)$, 
which counts the zeros of odd order on the critical line, these exceptions are to be expected.
\end{rem}

\begin{rem}
The main step in the proof of the theorem is the application of Daboussi and Delange's result
(see Lemma \ref{lem2} below). It must be pointed out that, although not in $\mathcal S$, in the case of the Epstein zeta-function associated to certain positive definite binary
 quadratic forms and in the case of the ideal class zeta-functions associated to certain
 quadratic number fields, the crucial exponential sum estimate \eqref{est1e} was derived from
 the arithmetic information contained in their Dirichlet coefficients and thereby the
 infinitude of zeros on the critical line was established (see \cite{P-T} and \cite{KC-NR}
 respectively).
\end{rem}

The plan of this note is as follows.
In section \ref{asym} we shall derive an asymptotic expression for $\Delta_F (s)$
and state some growth estimates. Section \ref{estimates} will deal with some exponential lemmas
and a lower bound estimate for the first power mean of a Dirichlet series on the critical line.
In section \ref{proof}, we shall prove the theorem.

\section{\emph{Asymptotic expansion and order estimates}}\label{asym}

The well-known Stirling's formula for the $\Gamma$-function states
that in any fixed vertical strip 
$-\infty < \alpha \leq \sigma \leq \beta < \infty$,
$$\Gamma(\sigma+it)=(2\pi)^{1/2} t^{\sigma+it-1/2}e^{-\frac{\pi}{2} t + 
\frac{\pi}{2} i (\sigma-1/2)-i t} ( 1 + O(1/t) ) \quad \hbox{as} \ t \rightarrow \infty .
$$
Using this formula along with the  functional equation for $F(s)$, we obtain
\begin{equation}\label{asym1e}
\Delta_F(s)=\omega_1 (Q_1 t^{d/2})^{1-2\sigma-2it}
t^{-i \Im \xi}e^{i d t}(1+O(1/t)),
\end{equation}
where
$\omega_1=\omega e^{-\frac{\pi}{2}i(d/2+\Re \xi)
+i \Im \xi}\prod_{j=1}^r \lambda_j^{-2i \Im \mu_j},$
with $|\omega_1|=1$,
and $Q_1=Q \prod_{j=1}^r \lambda_j^{\lambda_j}.$
Note that $Q_1=\sqrt{q}\, (2 \pi)^{-d/2}$. 

\medskip

As $d=2$ and $\Im \xi=0$ in the present case, we get
\begin{equation}\label{asym2e}
\Delta_F(s) = \omega_1^{\prime} (Q_1 t)^{1-2\sigma-2it} e^{2it} (1+O(1/t)).
\end{equation}
where 
$ \omega_1^{\prime} =\omega e^{-\frac{\pi}{2}i (1+ \xi )
}\prod_{j=1}^r \lambda_j^{-2i \Im \mu_j}.$

\medskip

Thus, we have 
\begin{equation}\label{asym3e}
\Delta_F^{-1/2}(s)=\omega_2 (Q_1 t)^{\sigma-1/2+it}e^{-i t} (1+O(1/t)),
\end{equation}

where $\omega_2 = (\omega_1^{\prime})^{-1/2}$ with $|\omega_2|=1$.
\medskip
%\section{\emph{Order estimates}}\label{order}

Next, we use the following uniform convexity estimates which are easy to verify.

$$
F(s)  = \left\{ \begin{array}{ll}
 O(Q_1^{1-\sigma}t^{1-\sigma + \epsilon}) \quad 0 \le \sigma \le 1 \\
 O(t^{\epsilon}), \  \sigma>1.
 \end{array} \right.
$$

\noindent
From \eqref{asym3e} we get

\begin{equation}\label{ord1e}
\Delta_F^{-1/2}(s)F(s)= \left\{ \begin{array}{ll}
O(Q_1^{1/2}t^{1/2+\epsilon}) & 0 \le \sigma \le 1 \\
O(Q_1^{1/2} t^{\sigma - 1/2 + \epsilon}) & 1 < \sigma \le 1+\delta
\end{array} \right. 
\end{equation}

\medskip

\section{\emph{Estimates on oscillatory integrals}}\label{estimates}

We use the following result from Potter and Titchmarsh
\cite{P-T} on exponential integrals.
\begin{lem}\label{lem1}
Let
$$ J = \int_T^{T'} t^{\alpha} \left(\frac{t}{e \beta}\right)^{it} dt,$$
where $\alpha,\beta>0$, and $0<T\le T' \le 2T$. Then
\begin{eqnarray}
 J & = & O(T^{\alpha}/\log(T/\beta)) \quad \textrm{if} \quad \beta<T,\label{exp1e}\\
 J & = & (2\pi)^{1/2}\beta^{\alpha+1/2}e^{i\pi/4-i\beta}+O(T^{\alpha+2/5})+
O(T^{\alpha}/\log(\beta/T))+{} \label{exp2e}\\
 & & {}+O(T^{\alpha}/\log(T'/\beta)) \quad \textrm{if} \quad T<\beta<T',\nonumber\\
 J & = & O(T^{\alpha}/\log(\beta/T')) \quad \textrm{if} \quad T'<\beta,\label{exp3e}\\
 J & = & O(T^{\alpha+1/2}) \quad \textrm{in any case}.\label{exp4e}
\end{eqnarray}
\end{lem}

We will also use the following result on exponential sums from Daboussi
and Delange \cite{D-D}.
\begin{lem}\label{lem2}
Let f be a multiplicative arithmetical function satisfying the condition
$\sum_{n\le x}|f(n)|^2 = O(x)$. Then, for every irrational $\alpha$, we have
\begin{equation}\label{est1e}
\sum_{n \le x} f(n)e^{2\pi i n \alpha}=o(x).
\end{equation}
\end{lem}

\begin{rem}
We can show that the method of Daboussi and Delange in \cite{D-D} can be worked out more carefully to replace $o(x)$ 
 in \eqref{est1e} with $O(x(\log \log x)^{-1/2})$. In fact, under the additional condition
that $|f(p)| \le A$ for all primes $p$, Montgomery and Vaughan \cite{M-V} improved the 
$O$-estimate to $O(x \log^{-1} x)$.
\end{rem}

\begin{rem}
For $\alpha$ rational, the estimate \eqref{est1e} was established by Daboussi and Delange \cite{D-D} under the restriction that $|f(n)| \le 1$ and some other conditions.
However, in the general case, the Dirichlet coefficients of functions in $\mathcal{S}$ do not
 satisfy this restriction.
\end{rem}

Now we state Theorem 3 of \cite{B}  as:

\begin{lem}\label{lem3}
Let $  B(s) = \sum_{n=1}^{\infty} b_n n^{-s} $ be any Dirichlet series satisfying the following conditions:
\begin{itemize}
\item[(i)] not all $b_n$'s are zero;
\item[(ii)] the function can be continued analytically in $ \sigma \geq a, \  |t| \geq t_0 $, and in this
region $ B(s) = O( (|t|+10)^A ).$
\end{itemize}
Then for every $\epsilon > 0, $ we have 
$$  \int_T^{T+H} | B(\sigma + it)| dt \gg H $$
for all $ H \geq (\log T)^{\epsilon}, \  T\geq T_0(\epsilon),$ and $\sigma > a $.
\end{lem}

\begin{rem} Ramachandra showed (see \cite{Kram}, Chapter II ) that the first power mean of a generalized Dirichlet series
satisfying certain conditions can not be too small. Lemma \ref{lem3} is a particular case of this general theorem, which
 is quite useful in  obtaining lower bounds of the type \eqref{pro1e}, 
even in short-intervals. 
\end{rem}
\section{\emph{Proof of the main theorem}}\label{proof}

We define the function
$$Z_F(t)=\Delta_F(1/2 + i t)^{-1/2}F(1/2 +i t),
$$
where $ Z_F(t)$ is the analogue of Hardy's function $Z(t)$ in the theory
of Riemann's zeta-function. The functional equation implies that $ Z_F(t)$
is real for real $t$. Thus the zeros of $F(s)$ on the critical line correspond 
to the real zeros of $ Z_F(t)$.

\medskip

Suppose now that $ Z_F(t)$ has no zeros in the interval $ [T,2T]  $ where $  T \geq T_1 > 0,  \ T_1 $ 
is a sufficiently large absolute constant.

\medskip

Consider the integral
$$ I =  \int_{T}^{2T} Z_F(t) dt.
$$
We are going to estimate $ |I| $ from below and above to derive eventually
a contradiction. Since the integrand is of constant sign by our assumption, we have
$$\ |I| = \int_{T}^{2T} |Z_F(t)| dt. 
$$

\medskip

The lower bound
\begin{equation}\label{pro1e} |I| \gg T 
\end{equation}
follows on taking $ H = T $ in Lemma \ref{lem3}.

\medskip

As for the upper bound estimation, we first write the integral $I$
as 
\begin{equation}\label{pro2e}
 I = -i \int_{1/2+iT}^{1/2+i 2 T}\Delta_F^{-1/2}(s)F(s)ds.  
\end{equation}
Next, we move the line of integration to $\sigma = 1+ \delta, ( \delta > 0 $ 
is a small positive constant less than $1/10 $) and 
apply Cauchy's theorem to the integral 
$$
\int \Delta_F^{-1/2}(s)F(s)ds,
$$ along the rectangle with sides 
$\sigma=1/2$, $\sigma=1+\delta$, $t=T$ and $t=2 T$.

\medskip

By Cauchy's theorem and  the estimate \eqref{ord1e}, the integral  \eqref{pro2e}
reduces to
\begin{equation}\label{pro3e}
\int_T^{2T}\Delta_F^{-1/2}(1+\delta+it)F(1+\delta+it)dt,
\end{equation}
with an error $O( T^{1/2+\delta + \epsilon })=o(T)$ coming from the horizontal lines.

\medskip

Substituting the asymptotic formula \eqref{asym3e} in \eqref{pro3e}, we see that the expression \eqref{pro3e} is a
constant multiple of
\begin{equation}\label{pro4e}
 \sum_{n=1}^{\infty}\frac{a(n)}{n^{1+\delta}}
\int_T^{2T}(Q_1 t)^{1/2+\delta}e^{it\log(\frac{Q_1 t}{e n})}\big(1+O(1/t)\big) dt.
\end{equation}

\medskip

The $O-$term in \eqref{pro4e} is trivially $O(T^{1/2+\delta})=o(T)$.

\medskip

The problem, therefore, reduces to the estimation of the sum
\begin{equation}\label{pro5e}
Q_1^{1/2+\delta} \sum_{n=1}^{\infty}\frac{a(n)}{n^{1+\delta}}\int_T^{2T}t^{1/2+\delta}
\left(\frac{Q_1 t}{e n}\right)^{i t}dt.
\end{equation}

\medskip

The integral in \eqref{pro5e} is evaluated by exponential integral techniques (Lemma \ref{lem1}).

\medskip

First we subdivide the sum in \eqref{pro5e} into sub-intervals
$$1\le n \le Q_1 T-1,\qquad Q_1 T +1 \le n \le 2 Q_1 T-1,\qquad 2 Q_1 T+1\le n,
$$
and denote the corresponding sums over these ranges as 
$\Sigma_1,\Sigma_2,\Sigma_3$ respectively.

\medskip
Notice that there are atmost 4 integers $n$ from the range of \eqref{pro5e} which have
not been included in the above ranges. As $n \asymp Q_1 T$ 
(by which we mean $ Q_1 T \ll n \ll Q_1 T$)
for these integers and $a(n)=O(n^{\epsilon})$, their contribution to \eqref{pro5e} is 
$$O\bigg( \frac{n^{\epsilon}}{n^{1+\delta}}T^{1/2+\delta}T\bigg)=
O\big( T^{1/2+\epsilon}\big)=o(T).
$$

The sum $\Sigma_1$ is estimated by using \eqref{exp1e} (with $T'=2T$) to be
\begin{equation}\label{pro6e}
\Sigma_1=O\bigg(T^{1/2+\delta}\sum \frac{| a(n) | }{n^{1+\delta}\log(Q_1 T/n)}\bigg).
\end{equation}

We then further sub-divide the range of \eqref{pro6e} into two sub-sums as follows
$$  1 \le n \le Q_1 T/2, \qquad Q_1 T/2 < n \le Q_1 T-1,
$$ and denote the corresponding sums by
 $\Sigma_{11}$ and $\Sigma_{12}$ respectively.

\medskip

For estimating $\Sigma_{11}$, we note that $\log(Q_1 T/n)\ge\log 2$ in this range 
and hence we get 
$$ \Sigma_{11} = O\bigg( T^{1/2+\delta}\sum \frac{| a(n) | }{n^{1+\delta}}\bigg) = O\left( T^{1/2+\delta}\right)=o(T).
$$

\medskip

For estimating $\Sigma_{12}$, we use the inequality $\log(Q_1 T/n)\ge (Q_1 T-n)/Q_1 T$,
$n\asymp Q_1 T$ and $a(n)=O(n^{\epsilon})$ and obtain
\begin{equation}\label{pro7e}
 \Sigma_{12} = O\bigg( T^{1/2+\delta}\sum \frac{| a(n) | }{n^{1+\delta}}
\frac{Q_1 T}{Q_1 T-n}\bigg) = O\bigg( T^{1/2+\epsilon}
\sum \frac{1}{Q_1 T-n}\bigg).
\end{equation}
Observe that the last sum in \eqref{pro7e} is
$$\sum (Q_1 T-n)^{-1} = \sum_{1 \le k < Q_1 T/2}(k+f)^{-1}=O(\log T),
$$
where $f = Q_1 T-\lfloor Q_1 T\rfloor \ge 0$ is the fractional part of $Q_1 T$.
 Using this in \eqref{pro7e} gives
$\Sigma_{12} = O( T^{1/2+\epsilon}\log T) =o(T)$.
\medskip

Estimating $\Sigma_3$ is similar to that of $\Sigma_1$ except that we use \eqref{exp3e} instead of \eqref{exp1e}. Hence, we find in a similar manner that $\Sigma_3 = o(T)$.
 
\medskip

The sum  $\Sigma_2$ is estimated by using \eqref{exp2e} (with $T'=2T$) to be
\begin{eqnarray}\label{pro8e}
\Sigma_2 &=& C\sum a(n)e^{-in/Q_1}+
O\left( \sum \frac{ | a(n) | }{n^{1+\delta}}T^{9/10+\delta}\right)+{}\\
 & & {}+O\left(T^{1/2+\delta}
 \sum \frac{|a(n)|}{n^{1+\delta}}E_n\right),\nonumber
\end{eqnarray}
where $C=(2\pi/Q_1)^{\frac{1}{2}}e^{i\frac{\pi}{4}}$ is a constant and
$E_n = \log^{-1}(n/(Q_1 T)) + \log^{-1}(2 Q_1 T/n)$.

\medskip

Using the fact that $Q_1=\sqrt{q}(2\pi)^{-1}$, we get that 
the main exponential sum in \eqref{pro8e} is a constant multiple of
$$\sum_{Q_1 T+1 \le n \le 2 Q_1 T-1}a(n)e^{-2\pi i n/ \sqrt{q}} \ .
$$
By the assumption \eqref{int1e} in Theorem 1, the irrationality of 
$\sqrt{q}$ and Lemma \ref{lem2}, we get this sum to be $o(T)$.

\medskip

The first $O$-term in \eqref{pro8e} is estimated by using $n \asymp Q_1 T$ and $a(n)=O(n^{\epsilon})$ to be $O(T^{9/10+\epsilon})=o(T)$.

\medskip

The estimation of the second $O$-term in \eqref{pro8e} is done in the same way as was done for $\Sigma_{12}$ to conclude that it is $O(( T)^{1/2+\epsilon}\log T) =o(T)$.

\medskip

Collecting all the estimates, we have
\begin{equation}\label{pro9e}
\int_{T}^{2T} Z_F(t) dt = o(T).
\end{equation}

Thus from \eqref{pro1e} and \eqref{pro9e}, we derive a contradiction. This
completes the proof of the theorem.

\medskip

\noindent
{\bf Acknowledgments:} The authors wish to thank Professor K. Soundararajan for pointing 
out \cite{M-V} to us. The authors also wish to thank Professors R. Balasubramanian, 
A. Ivi\'c and A. Perelli for their useful suggestions. Finally the authors express 
their gratitude to the referee for promptly reviewing the manuscript and reverting
back with some valuable suggestions.

\end{document}